\documentclass[11pt]{article}
\usepackage{}
\usepackage{mathrsfs}
\usepackage{amsfonts}

\usepackage{CJK}
\usepackage{amsfonts,amssymb,amsmath,mathrsfs}

\usepackage{color,latexsym,amsfonts}
 \newtheorem{theorem}{Theorem}[section]
 \newtheorem{lemma}[theorem]{Lemma}
 \newtheorem{corollary}[theorem]{Corollary}
 \newtheorem{proposition}[theorem]{Proposition}
 
 \newtheorem{Definition}[theorem]{Definition}
 \newtheorem{remark}[theorem]{Remark}
 \newtheorem{condition}[theorem]{Condition}

 \def\blemma{\begin{lemma}\sl{}\def\elemma{\end{lemma}}}
 \def\btheorem{\begin{theorem}\sl{}\def\etheorem{\end{theorem}}}

 \def\bremark{\begin{remark}\sl{}\def\eremark{\end{remark}}}

 \def\beqlb{\begin{eqnarray}}\def\eeqlb{\end{eqnarray}}
 \def\beqnn{\begin{eqnarray*}}\def\eeqnn{\end{eqnarray*}}

 \def\<{\langle}\def\>{\rangle}

 \def\eqref#1{{\rm(\ref{#1})}}

\def\D{\textup{D}}

\def\d{\textup{d}}

\def\e{\textup{e}}

\def\fin{\hfill$\square$}

\def\newdot{{\kern.8pt\cdot\kern.8pt}}

\def\R{\mathbb{R}}
\def\E{\mathbb{E}}
\def\P{\mathbb{P}}

\def\D{\mathbb{D}}

\def\<{\langle}
\def\>{\rangle}
\def\Proof.{\noindent{\bf Proof.}}

\begin{document}

\

\noindent{}

\bigskip\bigskip

\centerline{\Large\bf Derivative formulas and applications for degenerate SDEs}

\smallskip

\centerline{\Large\bf with fractional noises
\footnote{Supported by the National Natural Science Foundation of China (Grant No. 11501009, 11371029),
the Natural Science Foundation of Anhui Province (Grant No. 1508085QA03), and the Distinguished Young Scholars Foundation of Anhui Province (Grant No. 1608085J06).}
}

\smallskip
\
\bigskip\bigskip

\centerline{Xiliang Fan}

\bigskip

\centerline{Department of Statistics, Anhui Normal University, Wuhu 241003, China}
\smallskip
\centerline{fanxiliang0515@163.com}

\smallskip

\bigskip\bigskip

{\narrower{\narrower

\noindent{\bf Abstract.} For degenerate stochastic differential equations driven by fractional Brownian motions with Hurst parameter $H>1/2$,
the derivative formulas are established by using Malliavin calculus and coupling method, respectively.
Furthermore, we find some relation between these two approaches.
As applications, the (log) Harnack inequalities and the hyperbounded property are presented.}


\bigskip



\section{Introduction}

\setcounter{equation}{0}

The derivative formula initiated in \cite{Bismut84} is a powerful tool for stochastic analysis.
This, together with the integration by parts formula \cite{Driver97}, enables one to derive regular estimates on the commutator,
which is important for the study of flow properties \cite{Fang&Li&Luo11}.
Recently, by using a coupling technique or Malliavin calculus,
the derivative formulas have been extended and applied to various models for the study of L\'{e}vy processes.
One can see, for instance, \cite{Wang11a} for SDEs;
\cite{Bao&Wang&Yuan13b,Guillin&Wang12,Priola06,Wang&Zhang13,Zhang10a} for degenerate SDEs;
\cite{Bao&Wang&Yuan13c,Dong&Xie10,Elworthy&Li94,Wang&Xu10a,ZhangSQ13a} for SPDEs.
We remark that, the derivative formulas for damping stochastic Hamiltonian systems have been established
in \cite{Zhang10a} and \cite{Guillin&Wang12} respectively, where the degenerate part is linear.
Afterwards, Wang and Zhang \cite{Wang&Zhang13} extended the results derived in \cite{Zhang10a,Guillin&Wang12} to the non-linear degenerate case.
But, to the best of our knowledge, explicit derivative formula for degenerate SDEs in a non-L\'{e}vy context is not yet available.

In this paper, we are concerned with degenerate SDEs driven by fractional Brownian motions,
whose noises are not Le\'{v}y processes and even more not semimartingale when $H\neq1/2$.
Now, in essence there exist two mainly different integration theories to define and study the fractional stochastic integral:
the pathwise Riemann-Stieltjes integral originally due to \cite{Young36a} and the divergence operator (or the Skorohod integral)
(e.g. \cite{Alos&Mazet&Nualart01a,Decreusefond&Ustunel98a}) defined as the adjoint of the derivative operator in the framework of the Malliavin calculus.
Then, there are numerous works to investigate SDEs driven by fractional Brownian motions.
For instance, \cite{Coutin&Qian02a,Nualart&Rascanu02a} proved the existence and uniqueness result;
\cite{Baudoin&Ouyang11,Nourdin&Simon06,Nualart&Saussereau09}, \cite{Hairer05,Hairer&Pillai11},\cite{Saussereau12} and \cite{Baudoin&Ouyang&Tindel11a,Fan15b}
studied the distributional regularities, the ergodicity, the Talagrand transportation inequalities and the logarithmic Sobolev inequalities
for the solutions, respectively.
For other results involved with paths regularity of fractional Brownian motions,
one may refer to \cite{Xiao11a,Jaramillo&Nualart17a, Yan16a} and references therein.
Recently, in the previous papers \cite{Fan13b} and \cite{Fan13a},
we obtained derivative formulas for SDEs with fractional noises for $H<1/2$ and $H>1/2$, respectively.
Motivated by the work \cite{Wang&Zhang13}, where the derivative formulas were shown for the stochastic Hamiltonian system by using Malliavin calculus,
we will be able to establish the derivative formulas for degenerate SDEs with fractional noises ($H>1/2$),
which will imply Harnack type inequalities as well as hyperbounded property.
That is the main purpose of this article.

The rest of the paper is organized as follows.
In Section 2, we recall some basic results about fractional calculus and fractional Brownian motion.
The derivative formulas by means of Malliavin calculus and coupling argument will be addressed in Section 3 and 4, respectively.
In Section 5, with helps of these formulas,
we present some applications to the dimensional-free Harnack type inequalities and the hyperbounded property.

\section{Preliminaries}

\setcounter{equation}{0}

\subsection{Fractional Integrals and Derivatives}

\setcounter{equation}{0}

For later use, we introduce some basic facts about fractional calculus,
which can be found in \cite{Samko&Kilbas&Marichev}.

Let $a,b\in\R$ with $a<b$.
For $\alpha>0$ and $f\in L^1(a,b)$,
the left-sided (resp. right-sided) fractional Riemann-Liouville integral of $f$ of order $\alpha$ on $[a,b]$ is defined as
\beqnn
I_{a+}^\alpha f(x)=\frac{1}{\Gamma(\alpha)}\int_a^x\frac{f(y)}{(x-y)^{1-\alpha}}\d y,\ \
\left(\mbox{resp.}~I_{b-}^\alpha f(x)=\frac{(-1)^{-\alpha}}{\Gamma(\alpha)}\int_x^b\frac{f(y)}{(y-x)^{1-\alpha}}\d y\right),
\eeqnn
where $x\in(a,b)$ a.e., $(-1)^{-\alpha}=\e^{-i\alpha\pi},\Gamma$ stands for the Euler function.
In particular, when $\alpha=n\in\mathbb{N}$, they reduced to the usual $n$-order iterated integrals.

Let $p\geq1$.
If $f\in I_{a+}^\alpha(L^p)$ (resp. $I_{b-}^\alpha(L^p)$) and $0<\alpha<1$,
then the Weyl derivative reads as follow
\beqnn
D_{a+}^\alpha f(x)=\frac{1}{\Gamma(1-\alpha)}\left(\frac{f(x)}{(x-a)^\alpha}+\alpha\int_a^x\frac{f(x)-f(y)}{(x-y)^{\alpha+1}}\d y\right)\\
\left(\mbox{resp.}\ \ D_{b-}^\alpha f(x)=\frac{(-1)^\alpha}{\Gamma(1-\alpha)}\left(\frac{f(x)}{(b-x)^\alpha}+\alpha\int_x^b\frac{f(x)-f(y)}{(y-x)^{\alpha+1}}\d y\right)\right),
\eeqnn
where the convergence of the integrals at the singularity $y=x$ holds pointwise for almost all $x$ if $p=1$ and in $L^p$-sense if $1<p<\infty$.

Suppose that $f\in C^\lambda(a,b)$ (the set of $\lambda$-H\"{o}lder continuous functions on $[a,b]$) and $g\in C^\mu(a,b)$ with $\lambda+\mu>1$.
By \cite{Young36a}, the Riemann-Stieltjes integral $\int_a^bf\d g$ exists.
In \cite{Zahle98a}, Z\"{a}hle provides an explicit expression for the integral $\int_a^bf\d g$ in terms of fractional derivatives.
Let $\lambda>\alpha$ and $\mu>1-\alpha$.
Then the Riemann-Stieltjes integral can be expressed as
\beqlb\label{2.0}
\int_a^bf\d g=(-1)^\alpha\int_a^bD_{a+}^\alpha f(t)D_{b-}^{1-\alpha}g_{b-}(t)\d t,\nonumber
\eeqlb
where $g_{b-}(t)=g(t)-g(b)$.
This can be regarded as fractional integration by parts formula.

\subsection{Fractional Brownian Motion}

\setcounter{equation}{0}

In this part, we shall recall some important definitions and results concerning the fractional Brownian motion.
For a deeper discussion, we refer the reader to \cite{Alos&Mazet&Nualart01a,Biagini&Hu08a,Decreusefond&Ustunel98a} and \cite{Nualart06a}.

Let $H\in(1/2,1)$.
The $d$-dimensional fractional Brownian motion with Hurst parameter $H$ on the probability space $(\Omega,\mathscr{F},\mathbb{P})$
can be defined as the centered Gauss process $B^H=\{B_t^H, t\in[0,T]\}$ with covariance function $\E\left(B_t^{H,i}B_s^{H,j}\right)=R_H(t,s)\delta_{i,j}$,
where
\beqnn
 R_H(t,s)=\frac{1}{2}\left(t^{2H}+s^{2H}-|t-s|^{2H}\right).
\eeqnn

By the above covariance function, one can show that $\E|B_t^{H,i}-B_s^{H,i}|^p=C(p)|t-s|^{pH},\ \forall p\geq 1$.
Consequently, by the Kolmogorov continuity criterion
$B^{H,i}$ have $(H-\epsilon)$-order H\"{o}lder continuous paths for all $\epsilon>0,\ i=1,\cdot\cdot\cdot,d$.

For each $t\in[0,T]$, let $\mathcal {F}_t$ be the $\sigma$-algebra generated by the random variables $\{B_s^H:s\in[0,t]\}$ and the
$\mathbb{P}$-null sets.

Denote $\mathscr{E}$ by the set of step functions on $[0,T]$.
Let $\mathcal {H}$ be the Hilbert space defined as the closure of
$\mathscr{E}$ with respect to the scalar product
\beqnn
\langle (I_{[0,t_1]},\cdot\cdot\cdot,I_{[0,t_d]}),(I_{[0,s_1]},\cdot\cdot\cdot,I_{[0,s_d]})\rangle_\mathcal {H}=\sum\limits_{i=1}^dR_H(t_i,s_i).
\eeqnn
By bounded linear transform theorem,
the mapping $(I_{[0,t_1]},\cdot\cdot\cdot,I_{[0,t_d]})\mapsto\sum_{i=1}^dB_{t_i}^{H,i}$ can be extended to an isometry between $\mathcal {H}$ and the Gaussian space $\mathcal {H}_1$ associated with $B^H$.
Denote this isometry by $\phi\mapsto B^H(\phi)$.

On the other hand, from \cite{Decreusefond&Ustunel98a},
we know that the covariance kernel $R_H(t,s)$ has the following integral representation
\beqnn
 R_H(t,s)=\int_0^{t\wedge s}K_H(t,r)K_H(s,r)\d r,
\eeqnn
where $K_H$ is a square integrable kernel given by
\beqnn
K_H(t,s)=\Gamma\left(H+\frac{1}{2}\right)^{-1}(t-s)^{H-\frac{1}{2}}F\left(H-\frac{1}{2},\frac{1}{2}-H,H+\frac{1}{2},1-\frac{t}{s}\right),
\eeqnn
in which $F(\cdot,\cdot,\cdot,\cdot)$ is the Gauss hypergeometric function (for details see \cite{Decreusefond&Ustunel98a} or \cite{Nikiforov&Uvarov88}).

Now, define the linear operator $K_H^*:\mathscr{E}\rightarrow L^2([0,T],\R^d)$ by
\beqnn
(K_H^*\phi)(s)=K_H(T,s)\phi(s)+\int_s^T(\phi(r)-\phi(s))\frac{\partial K_H}{\partial r}(r,s)\d r.
\eeqnn
By integration by parts, it is easy to see that this can be rewritten as
\beqnn
(K_H^*\phi)(s)=\int_s^T\phi(r)\frac{\partial K_H}{\partial r}(r,s)\d r.
\eeqnn
Due to \cite{Alos&Mazet&Nualart01a}, for all $\phi,\psi\in\mathscr{E}$,
there holds $\langle K_H^*\phi,K_H^*\psi\rangle_{L^2([0,T],\R^d)}=\langle\phi,\psi\rangle_\mathcal {H}$
and then $K_H^*$ can be extended to an isometry between $\mathcal{H}$ and $L^2([0,T],\R^d)$.
Hence, according to \cite{Alos&Mazet&Nualart01a} again,
the process $\{W_t=B^H((K_H^*)^{-1}{\rm I}_{[0,t]}),t\in[0,T]\}$ is a Wiener process,
and $B^H$ has the following integral representation
\beqnn
 B^H_t=\int_0^tK_H(t,s)\d W_s.
\eeqnn

Define the operator $K_H: L^2([0,T],\mathbb{R}^d)\rightarrow I_{0+}^{H+1/2 }(L^2([0,T],\mathbb{R}^d))$ by
\beqnn
 (K_H f)(t)=\int_0^tK_H(t,s)f(s)\d s.
\eeqnn
By \cite{Decreusefond&Ustunel98a}, it is an isomorphism and for each $f\in L^2([0,T],\mathbb{R}^d)$,
\beqnn
 (K_H f)(s)=I_{0+}^{1}s^{H-1/2}I_{0+}^{H-1/2}s^{1/2-H}f.
\eeqnn
As a consequence, for every $h\in I_{0+}^{H+1/2}(L^2([0,T],\R^d))$,
the inverse operator $K_H^{-1}$ is of the  form
\beqlb\label{2.1}
(K_H^{-1}h)(s)=s^{H-1/2}D_{0+}^{H-1/2}s^{1/2-H}h'.
\eeqlb

The remaining part will be devoted to the Malliavin calculus of fractional Brownian motion.

Let $\Omega$ be the canonical probability space $C_0([0,T],\R^d)$, the set of continuous functions,
null at time $0$, equipped with the supremum norm.
Let $\P$ be the unique probability measure on $\Omega$
such that the canonical process $\{B^H_t; t\in[0,T]\}$ is a $d$-dimensional fractional Brownian motion with Hurst parameter $H$.
Then, the injection $R_H=K_H\circ K_H^*:\mathcal{H}\rightarrow\Omega$ embeds $\mathcal{H}$ densely into $\Omega$ and
$(\Omega,\mathcal{H},\P)$ is an abstract Wiener space in the sense of Gross.
In the sequel we will make this assumption on the underlying probability space.

Denote $\mathcal {S}$ by the set of smooth and cylindrical random variables of the form
\beqnn
F=f(B^H(\phi_1),\cdot\cdot\cdot,B^H(\phi_n)),
\eeqnn
where $n\geq 1, f\in C_b^\infty(\mathbb{R}^n)$, the set of $f$ and all its partial derivatives are bounded, $\phi_i\in\mathcal{H}, 1\leq i\leq n$.
The Malliavin derivative of $F$, denoted by $\mathbb{D}F$, is defined as the $\mathcal {H}$-valued random variable
\beqnn
\mathbb{D}F=\sum_{i=1}^n\frac{\partial f}{\partial x_i}(B^H(\phi_1),\cdot\cdot\cdot,B^H(\phi_n))\phi_i.
\eeqnn
For any $p\geq 1$, we define the Sobolev space $\mathbb{D}^{1,p}$ as the completion of $\mathcal {S}$ with respect to the norm
\beqnn
\|F\|_{1,p}^p=\mathbb{E}|F|^p+\mathbb{E}\|\mathbb{D}F\|^p_{\mathcal {H}}.
\eeqnn
While we will denote by $\delta$ and $\mathrm{Dom}\delta$ the divergence operator of $\mathbb{D}$ and its domain.
Let us finish this section by giving a transfer principle that connects the divergence operators of both processes $B^H$ and $W$.

\btheorem\label{Theorem2.1}\cite[Proposition 5.2.2]{Nualart06a}
$\mathrm{Dom}\delta=(K_H^*)^{-1}(\mathrm{Dom}\delta_W)$,
and for any $\mathcal{H}$-valued random variable $u$ in $\mathrm{Dom}\delta$ we have $\delta(u)=\delta_W(K_H^*u)$,
where $\delta_W$ denotes the divergence operator with respect to $W$.
\etheorem

\bremark\label{Remark2.1}
The above theorem, together with \cite[Proposition 1.3.11]{Nualart06a}, yields that
if $K_H^*u\in L_a^2([0,T]\times\Omega,\R^d)$ (the closed subspace of $L^2([0,T]\times\Omega,\R^d)$ formed by the adapted processes),
then $u\in\mathrm{Dom}\delta$.
\eremark

\section{Derivative Formulas by Malliavin Calculus}

\setcounter{equation}{0}

The objective of this section is to study the following degenerate SDE with fractional noise
\begin{equation}\label{3.1}
\begin{cases}
 \textnormal\d X_t=b_t(X_t,Y_t)\d t,\ \ \ \ \ \ \ \ \ \ \ \ \ \ \ \ \ X_0=x\in\R^{d_1},\\
 \textnormal\d Y_t=\tilde{b}_t(X_t,Y_t)\d t+\sigma(t)\d B_t^H, \ \ ~Y_0=y\in\R^{d_2},
\end{cases}
\end{equation}
where $b:[0,T]\times\R^{d_1}\times\R^{d_2}\rightarrow\R^{d_1}, \tilde{b}:[0,T]\times\R^{d_1}\times\R^{d_2}\rightarrow\R^{d_2},
\sigma:[0,T]\rightarrow\R^{d_2}\times\R^d$ and $H>1/2$.

Remark that, as in the Brownian motion case, when taking the special choices of $b,\tilde{b}$ and $\sigma$,
this model will be reduced to stochastic Hamiltonian system with fractional noises.
We shall use $(X_t^z,Y_t^z)$ to denote the solution with the initial value $z:=(x,y)\in\R^{d_1+d_2}:=\R^{d_1}\times\R^{d_2}$.
The associated operator $P_t$ is defined by
$$P_tf(z)=\E f(X_t^z,Y_t^z),\ t>0, z\in\R^{d_1+d_2},f\in\mathcal{B}_b(\R^{d_1+d_2}),$$
where $\mathcal{B}_b(\R^{d_1+d_2})$ is the set of all bounded measurable functions on $\R^{d_1+d_2}.$
Besides, set $Z_t^z=(X_t^z,Y_t^z), 0\leq t\leq T$ and for $f\in C^\alpha([0,T];\R^m)$,
put
\beqnn
\|f\|_\alpha:=\sup\limits_{s\neq t,s,t\in[0,T]}\frac{|f(t)-f(s)|}{|t-s|^\alpha}.
\eeqnn

We begin with the assumption (H1)
\begin{itemize}
\item[(i)] $b$ is Lipschitz continuous in $z$:
$$|(b_t,\tilde{b}_t)(z_1)-(b_t,\tilde{b}_t)(z_2)|\leq L|z_1-z_2|, \ \forall t\in[0,T], z_1,z_2\in\R^{d_1+d_2},$$
and for each $z\in\R^{d_1+d_2}, (b_\cdot,\tilde{b}_\cdot)(z)$ is Lipschitz continuous;
\par

\item[(ii)] $\sigma$ is H\"{o}lder continuous of order $\delta\in\left((1-H)\vee(H-1/2),1\right]$:
$$|\sigma(t)-\sigma(s)|\leq\tilde{L}|t-s|^\delta, \ \forall t,s\in[0,T],$$
and $\sigma\sigma^*$ is invertible so that $(\sigma\sigma^*)^{-1}$ is bounded;

where $L$ and $\tilde{L}$ are both nonnegative constants.
\end{itemize}

Due to \cite[Theorem 2.1]{Nualart&Rascanu02a},
the condition (H1) ensures that there exists a unique adapted solution to the equation \eqref{3.1}.
We aim to establish the Bismut type derivative formulas for the operator $P_T$ by means of Malliavin calculus.
That is, for any $v=(v_1,v_2)\in\R^{d_1+d_2}$,
we are to find $h\in\mathrm{Dom}\delta$ such that
\beqlb\label{Rev-add3.3}
\nabla_v P_T f(z)=\E(f(Z^z_T)\delta(h)),\ f\in C_b^1(\R^{d_1+d_2})
\eeqlb
To this end, for a stochastic process with differentiable paths $(\tilde{g}(t))_{0\leq t\leq T}$,
let $g$ solve the following linear equation
\beqlb\label{add3.4}
g(t)=v_1+\int_0^t\nabla^1b_s(Z^z_s)g(s)\d s+\int_0^t\nabla^2b_s(Z^z_s)\tilde{g}(s)\d s,
\eeqlb
and then put
\beqlb\label{add3.1}
(R_H h)(t)
=\int_0^t\sigma^*(s)((\sigma\sigma^*)(s))^{-1}\left(\nabla^1\tilde{b}_s(Z^z_s)g(s)+\nabla^2\tilde{b}_s(Z^z_s)\tilde{g}(s)-\tilde{g}'(s)\right)\d s.
\eeqlb
Here we use the notations $\nabla^1$ and $\nabla^2$ to represent the gradient operators on $\R^{d_1}$ and $\R^{d_2}$ respectively,
and recall that the operator $R_H$ is defined as $K_H\circ K_H^*$ in Section 2.
Our main result reads as follow.

\btheorem\label{Th3.1}
Assume (H1).
For $v=(v_1,v_2)\in\R^{d_1+d_2}$,
let $(\tilde{g}(t))_{0\leq t\leq T}$ be a stochastic process with differentiable paths such that $\tilde{g}(0)=v_2, \tilde{g}(T)=0$,
and $g,h$ given in \eqref{add3.4},\eqref{add3.1}, respectively.
If $h\in\mathrm{Dom}\delta$ and $g(T)=0$,
then there holds \eqref{Rev-add3.3}.
Furthermore, if $K_H^*h\in L_a^2([0,T]\times\Omega,\R^d)$, then we have
\beqlb\label{add3.0}
&&\delta(h)=\delta_W(K_H^*h)\nonumber\\
&&=\int_0^T\langle K^*_Hh(t),\d W(t)\rangle\nonumber\\
&&=\int_0^T\left\langle K_H^{-1}\left(\int_0^\cdot\sigma^*(s)((\sigma\sigma^*)(s))^{-1}\left(\nabla^1\tilde{b}_s(Z^z_s)g(s)+\nabla^2\tilde{b}_s(Z^z_s)\tilde{g}(s)-\tilde{g}'(s)\right)\right)(t),\d W(t)\right\rangle.\nonumber\\
\eeqlb
\etheorem

\emph{Proof.}
By \cite[Lemma 3.1 and the proof of Proposition 3.3]{Fan13a} we have $X^{z,i}_t,Y^{z,j}_t\in\D^{1,2}, 1\leq i\leq d_1, 1\leq j\leq d_2$,
and moreover their Malliavin derivatives solve the following equation:
for each $h\in\mathcal {H}$,
\begin{equation}\label{3.2}
\left\{
\begin{array}{ll}
\langle\mathbb{D}X^{z,i}_t,h\rangle_\mathcal {H}=\sum\limits_{k=1}^{d_1}\int_0^t(\nabla b_s(Z^z_s))_{ik}
    \langle\mathbb{D}X^{z,k}_s,h\rangle_\mathcal {H}\d s+
    \sum\limits_{k=1}^{d_2}\int_0^t(\nabla b_s(Z^z_s))_{ik}\langle\mathbb{D}Y^{z,k}_s,h\rangle_\mathcal {H}\d s,\\
\langle\mathbb{D}Y^{z,j}_t,h\rangle_\mathcal {H}=\sum\limits_{k=1}^{d_1}\int_0^t(\nabla \tilde{b}_s(Z^z_s))_{jk}
    \langle\mathbb{D}X^{z,k}_s,h\rangle_\mathcal {H}\d s+
     \sum\limits_{k=1}^{d_2}\int_0^t(\nabla \tilde{b}_s(Z^z_s))_{jk}\langle\mathbb{D}Y^{z,k}_s,h\rangle_\mathcal {H}\d s\\
     \ \ \ \ \ \ \ \ \ \ \ \ \ \ \ \ \ \ \ +\sum\limits_{k=1}^d\int_0^t\sigma_{jk}(s)\d (R_H h)^k(s).
\end{array} \right.
\end{equation}
Again by the arguments of \cite[Lemma 3.1 and Proposition 3.3]{Fan13a}, we conclude that
$\mathbb{D}_{R_Hh}Z^z_t=\frac{\d}{\d\epsilon}\Big|_{\epsilon=0}Z^z_t(w+\epsilon R_Hh),\ h\in\mathcal {H},$ satisfies
\begin{equation}\label{3.3}
\mathbb{D}_{R_Hh}Z^z_t=\langle\mathbb{D}Z^z_t,h\rangle_\mathcal{H}.
\end{equation}
Hence, we can reformulate \eqref{3.2} as
\begin{equation}\label{3.4}
\left\{
\begin{array}{ll}
\mathbb{D}_{R_Hh}X^z_t=\int_0^t(\nabla^1b_s)(Z^z_s)\mathbb{D}_{R_Hh}X^z_s\d s
   +\int_0^t(\nabla^2b_s)(Z^z_s)\mathbb{D}_{R_Hh}Y^z_s\d s,\\
\mathbb{D}_{R_Hh}Y^z_t=\int_0^t(\nabla^1\tilde{b}_s)(Z^z_s)\mathbb{D}_{R_Hh}X^z_s\d s
+\int_0^t(\nabla^2\tilde{b}_s)(Z^z_s)\mathbb{D}_{R_Hh}Y^z_s\d s+\int_0^t\sigma(s)\d (R_H h)(s).
\end{array} \right.
\end{equation}

On the other hand, by a direct calculation we show that the directional derivative processes satisfy that,
for any $v_1\in\R^{d_1}$ and $v_2\in\R^{d_2}$,
\begin{equation}\label{3.5}
\left\{
\begin{array}{ll}
\nabla_{v_1}^1X^z_t=v_1+\int_0^t\nabla^1b_s(Z^z_s)\nabla_{v_1}^1X^z_s\d s+\int_0^t\nabla^2b_s(Z^z_s)\nabla_{v_1}^1Y^z_s\d s,\\
\nabla_{v_1}^1Y^z_t=\int_0^t\nabla^1\tilde{b}_s(Z^z_s)\nabla_{v_1}^1X^z_s\d s+\int_0^t\nabla^2\tilde{b}_s(Z^z_s)\nabla_{v_1}^1Y^z_s\d s,
\end{array} \right.
\end{equation}
and
\begin{equation}\label{3.6}
\left\{
\begin{array}{ll}
\nabla_{v_2}^2X^z_t=\int_0^t\nabla^1b_s(Z^z_s)\nabla_{v_2}^2X^z_s\d s+\int_0^t\nabla^2b_s(Z^z_s)\nabla_{v_2}^2Y^z_s\d s,\\
\nabla_{v_2}^2Y^z_t=v_2+\int_0^t\nabla^1\tilde{b}_s(Z^z_s)\nabla_{v_2}^2X^z_s\d s+\int_0^t\nabla^2\tilde{b}_s(Z^z_s)\nabla_{v_2}^2Y^z_s\d s.
\end{array} \right.
\end{equation}
By $\nabla_vX^z_t=\nabla^1_{v_1}X^z_t+\nabla^2_{v_2}X^z_t, \nabla_vY^z_t=\nabla^1_{v_1}Y^z_t+\nabla^2_{v_2}Y^z_t$,
combining \eqref{3.5} with \eqref{3.6} yields that
\begin{equation}\label{3.7}
\left\{
\begin{array}{ll}
\nabla_vX^z_t=v_1+\int_0^t\nabla^1b_s(Z^z_s)\nabla_vX^z_s\d s+\int_0^t\nabla^2b_s(Z^z_s)\nabla_vY^z_s\d s,\\
\nabla_vY^z_t=v_2+\int_0^t\nabla^1\tilde{b}_s(Z^z_s)\nabla_vX^z_s\d s+\int_0^t\nabla^2\tilde{b}_s(Z^z_s)\nabla_vY^z_s\d s.
\end{array} \right.
\end{equation}

Note that, taking into account $\tilde{g}(0)=v_2, \eqref{add3.4}$ and $\eqref{add3.1}$, we obtain
\begin{equation}\label{3.8}
\left\{
\begin{array}{ll}
g(t)=v_1+\int_0^t\nabla^1b_s(Z^z_s)g(s)\d s+\int_0^t\nabla^2b_s(Z^z_s)\tilde{g}(s)\d s,\\
\tilde{g}(t)=v_2+\int_0^t\nabla^1\tilde{b}_s(Z^z_s)g(s)\d s+\int_0^t\nabla^2\tilde{b}_s(Z^z_s)\tilde{g}(s)\d s-\int_0^t\sigma(s)\d R_Hh(s).
\end{array} \right.
\end{equation}
This, together with \eqref{3.4}, leads to
\begin{equation}\label{3.9}
\left\{
\begin{array}{ll}
\mathbb{D}_{R_Hh}X^z_t+g(t)=v_1+\int_0^t(\nabla^1b_s)(Z^z_s)(\mathbb{D}_{R_Hh}X^z_s+g(s))\d s
   +\int_0^t(\nabla^2b_s)(Z^z_s)(\mathbb{D}_{R_Hh}Y^z_s+\tilde{g}(s))\d s,\\
\mathbb{D}_{R_Hh}Y^z_t+\tilde{g}(t)=v_2+\int_0^t(\nabla^1\tilde{b}_s)(Z^z_s)(\mathbb{D}_{R_Hh}X^z_s+g(s))\d s
+\int_0^t(\nabla^2\tilde{b}_s)(Z^z_s)(\mathbb{D}_{R_Hh}Y^z_s+\tilde{g}(s))\d s.
\end{array} \right.
\end{equation}
By \eqref{3.7}, \eqref{3.9} and the uniqueness of solutions of the ODE, we get
$$(\mathbb{D}_{R_Hh}X^z_t+g(t),\mathbb{D}_{R_Hh}Y^z_t+\tilde{g}(t))=(\nabla_vX^z_t,\nabla_vY^z_t),\ \forall t\in[0,T].$$
In particular, due to $g(T)=\tilde{g}(T)=0$, we have
$$(\mathbb{D}_{R_Hh}X^z_T,\mathbb{D}_{R_Hh}Y^z_T)=(\nabla_vX^z_T,\nabla_vY^z_T).$$
That is, $\mathbb{D}_{R_Hh}Z^z_T=\nabla_vZ^z_T.$
Hence, by $h\in\mathrm{Dom}\delta$, it follows that
\beqnn
&&\nabla_v P_T f(z)=\E(\nabla f(Z^z_T)\nabla_vZ^z_T)=\E(\nabla f(Z^z_T)\mathbb{D}_{R_Hh}Z^z_T)\\
&&=\E\mathbb{D}_{R_Hh}f(Z^z_T)=\E\langle\mathbb{D}f(Z^z_T),h\rangle_\mathcal {H}
=\E(f(Z^z_T)\delta(h)).
\eeqnn
Finally, if $K_H^*h\in L_a^2([0,T]\times\Omega,\R^d)$, then it follows by Remark \ref{Remark2.1} that $h\in\mathrm{Dom}\delta$.
Therefore, Theorem \ref{Theorem2.1} and \eqref{add3.1} allow to conclude the equation \eqref{add3.0}.
\fin

Next, we intend to apply Theorem \ref{Th3.1} with concrete choices of $\tilde{g}$ and then $g,h$.
To this end,
in the rest of the paper, we consider a special case of equation \eqref{3.1},
when $b(x,y)=Ax+By$.
That is, we consider the following fractional degenerate SDE
\begin{equation}\label{S3.1}
\begin{cases}
 \textnormal\d X_t=(AX_t+BY_t)\d t,\ \ \ \ \ \ \ \ \ \ \ \ \ \ \ \ \ X_0=x\in\R^{d_1},\\
 \textnormal\d Y_t=\tilde{b}_t(X_t,Y_t)\d t+\sigma(t)\d B_t^H, \ \ \ \ \ \ ~~Y_0=y\in\R^{d_2}.
\end{cases}
\end{equation}

For the equation \eqref{S3.1}, we impose additional condition on the coefficient $\tilde{b}$ needed to state our result: (H2)

$\nabla\tilde{b}$ is bounded and H\"{o}lder continuous:
$$|\nabla\tilde{b}_t(z_1)-\nabla\tilde{b}_t(z_2)|+|\nabla\tilde{b}_t(z_1)-\nabla\tilde{b}_s(z_1)|\leq K\left(|z_1-z_2|^\gamma+|t-s|^\varrho\right), \ \forall z_1,z_2\in\R^{d_1+d_2}, s,t\in[0,T],$$
where $\gamma\in(1-1/(2H),1],\varrho\in(H-1/2,1]$ and $K$ is a nonnegative constant.

We now have
\btheorem\label{Th3.2}
Assume (H1)(ii) and (H2).
Let $\alpha_i,i=1,2$ and $\rho$ be three functions such that
$\alpha_i\in C^2([0,T]),i=1,2$ and $\rho\in C([0,T])$ with $\alpha_1(0)=1,\alpha_1(T)=0, \alpha_2(0)=\alpha_2(T)=0$ and $\alpha_2(t)>0, \rho(t)>0, t\in(0,T]$ satisfying
\beqlb\label{Th3.2-0}
U_t:=\int_0^t\alpha_2(s)\e^{(T-s)A}BB^*\e^{(T-s)A^*}\d s\geq\rho(t)I_{d_1\times d_1}, \ t\in(0,T].
\eeqlb
Then, for $v=(v_1,v_2)\in\R^{d_1+d_2}$ and $ f\in C_b^1(\R^{d_1+d_2})$ there holds
\beqlb\label{Th3.2-1}
&&\nabla_v P_T f(z)\nonumber\\
&&=\E\left(f(Z^z_T)\int_0^T\left\langle K_H^{-1}\left(\int_0^\cdot\sigma^*(s)((\sigma\sigma^*)(s))^{-1}\left(\nabla^1\tilde{b}_s(Z^z_s)g(s)+\nabla^2\tilde{b}_s(Z^z_s)\tilde{g}(s)-\tilde{g}'(s)\right)\right)(t),\d W(t)\right\rangle\right),\nonumber\\
\eeqlb
where
\beqlb\label{Th3.2-2}
g(t)=\e^{tA}v_1+\int_0^t\e^{(t-s)A}B\tilde{g}(s)\d s,
\eeqlb
and
\beqlb\label{Th3.2-3}
\tilde{g}(t)&=&\alpha_1(t)v_2-\alpha_2(t)B^*\e^{(T-t)A^*}U_T^{-1}\e^{TA}v_1\nonumber\\
&&-\alpha_2(t)B^*\e^{(T-t)A^*}U_T^{-1}\int_0^T\alpha_1(s)\e^{(T-s)A}B v_2\d s.
\eeqlb
\etheorem

\bremark\label{addRemark(Th3.2)}
In order to ensure that \eqref{Th3.2-0} in the above theorem holds,
one needs to impose some non-degeneracy condition on the matrix $B$.
For instance, set the following Kalman rank condition:\\
\indent There exists an integer number $k_0\in[0,d_1-1]$ such that Rank$[B,AB,\cdots,A^{k_0}B]=d_1$.\\
When $K_0=0$, this condition reduces to Rank$[B]=d_1$.
Then \eqref{Th3.2-0} holds with $\alpha_2(t)=\frac{t(T-t)}{T^2}$
and $\rho(t)=\frac{C_1(t\wedge1)^{2(K_0+1)}}{T\e^{C_2T}}$,
where $C_i,i=1,2$ are both positive constants.
For more details, one can refer to \cite[Theorem 4.2]{Wang&Zhang13}.
\eremark

Before proving Theorem \ref{Th3.2},
we will first give a technical lemma concerning the estimation of the solution to \eqref{S3.1}.
The proof is identical to Lemma 3.4 proposed in \cite{Fan15a} and so we omit it.

\blemma\label{Le3.1}
Assume that (H1)(ii) and (H2) are fulfilled.
Then, there hold
\beqnn
\|Z^z\|_\infty\leq C(1+\|B^H\|_\lambda)
\eeqnn
and
\beqnn
|Z^z_t-Z^z_s|\leq C(|t-s|+\|B^H\|_\lambda|t-s|^\lambda), \ \forall s,t\in[0,T],
\eeqnn
\elemma
where $C$ is a positive constant and $\lambda$ is taken satisfying $\lambda\gamma>H-1/2$ and $\lambda\in(1-\delta,H)$.

\emph{Proof of Theorem \ref{Th3.2}.}
We observe first that $g(t)$ defined in \eqref{Th3.2-2} satisfies equation \eqref{add3.4} due to $\nabla^1b=A$ and $\nabla^2b=B$.
By \eqref{Th3.2-0}, it is easy to check that $U_t$ is invertible for all $t\in(0,T]$.
As a consequence, we are able to choose $\tilde{g}(t)$ as in \eqref{Th3.2-3} to obtain $g(T)=0$,
and then with the help of the values of $\alpha_i,i=1,2$ at the interval endpoints,
we deduce that $\tilde{g}(0)=v_2,\tilde{g}(T)=0$.
More precisely,
\beqnn
g(T)&=&\e^{TA}v_1+\int_0^T\e^{(T-s)A}B\tilde{g}(s)\d s\\
&=&\e^{TA}v_1+\int_0^T\alpha_1(s)\e^{(T-s)A}B v_2\d s
-\int_0^T\alpha_2(s)\e^{(T-s)A}BB^*\e^{(T-s)A^*}U_T^{-1}\e^{TA}v_1\d s\\
&&-\int_0^T\alpha_2(s)\e^{(T-s)A}BB^*\e^{(T-s)A^*}U_T^{-1}\int_0^T\alpha_1(s)\e^{(T-s)A}B v_2\d s\\
&=&0,
\eeqnn
and since $\alpha_1(0)=1, \alpha_1(T)=0$ and $\alpha_2(0)=0,\alpha_2(T)=0$,
by \eqref{Th3.2-3} we verify easily that $\tilde{g}(0)=v_2,\tilde{g}(T)=0$.

Now, to show \eqref{Th3.2-1}, by Theorem \ref{Th3.1} it remains to verify that $K_H^*h\in L_a^2([0,T]\times\Omega,\R^d)$, i.e.
$$K_H^{-1}\left(\int_0^\cdot\sigma^*(s)((\sigma\sigma^*)(s))^{-1}
\left(\nabla^1\tilde{b}_s(Z^z_s)g(s)+\nabla^2\tilde{b}_s(Z^z_s)\tilde{g}(s)-\tilde{g}'(s)\right)\right)\in L^2_a([0,T]\times\Omega,\R^d).$$
It is obvious that the operator $K^{-1}_H$ preserves the adaptability property.
Due to \eqref{2.1}, we have
\beqlb\label{Th3.2-pf1}
&&K_H^{-1}\left(\int_0^\cdot\sigma^*(s)((\sigma\sigma^*)(s))^{-1}
\left(\nabla^1\tilde{b}_s(Z^z_s)g(s)+\nabla^2\tilde{b}_s(Z^z_s)\tilde{g}(s)-\tilde{g}'(s)\right)\right)(t)\nonumber\\
&=:&K_H^{-1}\left(\int_0^\cdot\sigma^*(s)((\sigma\sigma^*)(s))^{-1}\eta(s)\right)(t)\nonumber\\
&=&t^{H-\frac{1}{2}}D^{H-\frac{1}{2}}_{0+}
\left[\cdot^{\frac{1}{2}-H}\sigma^*(\cdot)((\sigma\sigma^*)(\cdot))^{-1}
\eta(\cdot)\right](t)\nonumber\\
&=&\frac{1}{\Gamma(\frac{3}{2}-H)}
\Bigg[t^{\frac{1}{2}-H}\sigma^*(t)((\sigma\sigma^*)(t))^{-1}\eta(t)\nonumber\\
&&~~~~~~~~~~~~~~+\left(H-\frac{1}{2}\right)t^{H-\frac{1}{2}}\int_0^t\frac{t^{\frac{1}{2}-H}-s^{\frac{1}{2}-H}}{(t-s)^{\frac{1}{2}+H}}\sigma^*(s)((\sigma\sigma^*)(s))^{-1}\eta(s)\d s\nonumber\\
&&~~~~~~~~~~~~~~+\left(H-\frac{1}{2}\right)\int_0^t\frac{\sigma^*(t)((\sigma\sigma^*)(t))^{-1}-\sigma^*(s)((\sigma\sigma^*)(s))^{-1}}{(t-s)^{\frac{1}{2}+H}}\eta(t)\d s\nonumber\\
&&~~~~~~~~~~~~~~+\left(H-\frac{1}{2}\right)\int_0^t\sigma^*(s)((\sigma\sigma^*)(s))^{-1}\frac{\eta(t)-\eta(s)}{(t-s)^{\frac{1}{2}+H}}\d s\Bigg]\nonumber\\
&=:&\frac{1}{\Gamma(\frac{3}{2}-H)}[J_1(t)+J_2(t)+J_3(t)+J_4(t)].
\eeqlb
Note that by \eqref{Th3.2-0}, there holds $\|U_t^{-1}\|\leq1/\rho(t), t\in(0,T]$.
Therefore, combining this with \eqref{Th3.2-2}, \eqref{Th3.2-3} and the hypotheses,
we show that $\eta$ is bounded.
Moreover, noting that
\beqnn
\int_0^t\frac{s^{\frac{1}{2}-H}-t^{\frac{1}{2}-H}}{(t-s)^{\frac{1}{2}+H}}\d s
=t^{1-2H}\int_0^1\frac{r^{\frac{1}{2}-H}-1}{(1-r)^{\frac{1}{2}+H}}\d r<\infty
\eeqnn
and the fact: $\sigma^*(\sigma\sigma^*)^{-1}$ is $\delta$-order H\"{o}lder continuous,
we conclude that $J_i\in L^2([0,T]\times\Omega,\R^d), i=1,2,3.$
As for $\int_0^T|J_4(t)|^2\d t$, by a direct calculus it reduces to the following two terms:
\beqlb\label{Th3.2-pf2}
\int_0^T\left(\int_0^t\left|\frac{\nabla^i\tilde{b}_t(Z_t^z)-\nabla^i\tilde{b}_t(Z_s^z)}{(t-s)^{\frac{1}{2}+H}}\right|\d s\right)^2\d t,\ i=1,2.
\eeqlb
Hence, using (H2), Lemma \ref{Le3.1} and the Fernique theorem (see, for instance, \cite[Lemma 8]{Saussereau12}) along with \eqref{Th3.2-pf2},
we obtain $\int_0^T|J_4(t)|^2\d t<\infty$.
Thus the theorem is proved.
\fin

\section{Derivative Formulas by Coupling Method}

\setcounter{equation}{0}

In this section, our objective is to give derivative formula for \eqref{3.1} by coupling argument.
Recall that, in terms of Brownian motion, when $b$ is non-linear,
it seems very hard to construct successful couplings to derive derivative formulas and Harnack type inequalities.
In \cite{Wang&Zhang13}, the authors split $b$ into two terms: a linear term and a non-linear term,
and then controlled the latter part by the former part in some way.
In order to improve readability and to clarify the relations between Malliavin calculus and coupling method in dealing with derivative formulas,
we study the equation \eqref{S3.1}.
That is,
\begin{equation}\label{4.1}
\begin{cases}
 \textnormal\d X_t=(AX_t+BY_t)\d t,\ \ \ \ \ \ \ \ \ \ \ \ X_0=x\in\R^{d_1},\\
 \textnormal\d Y_t=\tilde{b}_t(X_t,Y_t)\d t+\sigma(t)\d B_t^H, \ \ ~Y_0=y\in\R^{d_2}.
\end{cases}
\end{equation}
For the convenience of the notations, let
$$\vartheta(t)=\sigma^*(t)((\sigma\sigma^*)(t))^{-1},\ \kappa(t)=\e^{tA}\tilde{x}+\int_0^t\e^{(t-s)A}B\tilde{h}(s)\d s, \ t\in[0,T],\tilde{x}\in\R^{d_1},$$
where $\tilde{h}$ is given \eqref{Th4.1-03} below.

Let us give now a statement for the derivative formula to the equation \eqref{4.1}.
We mention that the formula in the following theorem is indeed the same as Theorem \ref{Th3.2} (see Remark \ref{R4.2'} below),
but here, we use a different approach, i.e. coupling argument.

\btheorem\label{Th4.1-01}
Assume the same conditions as Theorem \ref{Th3.2} and let $\tilde{z}=(\tilde{x},\tilde{y})\in\R^{d_1+d_2}$.
Then
\beqlb\label{Th4.1-1}
\nabla_{\tilde{z}}P_T f(z)=\E\left[f(X_T^z,Y_T^z)M_T\right],\ \ f\in C_b^1(\R^{d_1+d_2}),
\eeqlb
where
\beqlb\label{Th4.1-02}
&&M_T\nonumber\\
&=&
\frac{H-\frac{1}{2}}{\Gamma(\frac{3}{2}-H)}\Bigg\{\frac{1}{H-\frac{1}{2}}\int_0^T\left\langle s^{\frac{1}{2}-H}
\vartheta(s)\left(\nabla\tilde{b}_s(X_s^z,Y_s^z)(\kappa(s),\tilde{h}(s))-\tilde{h}'(s)\right),\d W_s\right\rangle\nonumber\\
&&+\int_0^T\left\langle s^{H-\frac{1}{2}}
\vartheta(s)\left(\nabla\tilde{b}_s(X_s^z,Y_s^z)(\kappa(s),\tilde{h}(s))-\tilde{h}'(s)\right)
\int_0^s\frac{s^{\frac{1}{2}-H}-r^{\frac{1}{2}-H}}{(s-r)^{\frac{1}{2}+H}}\d r,\d W_s\right\rangle\nonumber\\
&&+\int_0^T\left\langle\left(\nabla\tilde{b}_s(X_s^z,Y_s^z)(\kappa(s),\tilde{h}(s))-\tilde{h}'(s)\right)
\int_0^s\left(\frac{s}{r}\right)^{H-\frac{1}{2}}
\frac{\vartheta(s)-\vartheta(r)}{(s-r)^{\frac{1}{2}+H}}\d r,\d W_s\right\rangle\nonumber\\
&&+\int_0^T\left\langle\int_0^s\left(\frac{s}{r}\right)^{H-\frac{1}{2}}
\vartheta(r)
\frac{\nabla\tilde{b}_s(X_s^z,Y_s^z)(\kappa(s),\tilde{h}(s))
-\nabla\tilde{b}_r(X_r^z,Y_r^z)(\kappa(r),\tilde{h}(r))}{(s-r)^{\frac{1}{2}+H}}\d r,\d W_s\right\rangle\nonumber\\
&&-\int_0^T\left\langle\int_0^s\left(\frac{s}{r}\right)^{H-\frac{1}{2}}
\vartheta(r)
\frac{\tilde{h}'(s)-\tilde{h}'(r)}{(s-r)^{\frac{1}{2}+H}}\d r,\d W_s\right\rangle\Bigg\}
\eeqlb
and
\beqlb\label{Th4.1-03}
\tilde{h}(t)=\alpha_1(t)\tilde{y}-\alpha_2(t)B^*\e^{(T-t)A^*}U_T^{-1}\e^{TA}\tilde{x}
-\alpha_2(t)B^*\e^{(T-t)A^*}U_T^{-1}\int_0^T\alpha_1(s)\e^{(T-s)A}B\tilde{y}\d s.
\eeqlb
\etheorem

\emph{Proof.}
For $\epsilon\in(0,1)$, let $(X_t^\epsilon,Y_t^\epsilon)$ solve the coupling equation
\begin{equation}\label{4.2}
\begin{cases}
 \textnormal\d X_t^\epsilon=(AX_t^\epsilon+BY_t^\epsilon)\d t,\ \ \ \ \ \ \ \ \ \ \ \ \ \ \ \ \ \ \ \ \ \ \ \ \ \ \ X_0^\epsilon=x+\epsilon\tilde{x}\in\R^{d_1},\\
 \textnormal\d Y_t^\epsilon=\tilde{b}_t(X_t^z,Y_t^z)\d t+\sigma(t)\d B_t^H+\epsilon\tilde{h}'(t)\d t, \ \ ~Y_0^\epsilon=y+\epsilon\tilde{y}\in\R^{d_2}.
\end{cases}
\end{equation}
Combining \eqref{4.1} with \eqref{4.2} yields that
\begin{equation}\label{4.3}
\begin{cases}
 \textnormal\d(X_t^\epsilon-X_t^z)=(A(X_t^\epsilon-X_t^z)+B(Y_t^\epsilon-Y_t^z))\d t,\ \ \ X_0^\epsilon-X_0^z=\epsilon\tilde{x}\in\R^{d_1},\\
 \textnormal\d(Y_t^\epsilon-Y_t^z)=\epsilon\tilde{h}'(t)\d t, \ \ \ \ \ \ \ \ \ \ \ \ \ \ \ \ \ \ \ \ \ \ \ \ \ \ \ \ \ \ \ \ \ \ Y_0^\epsilon-Y_0^z=\epsilon\tilde{y}\in\R^{d_2}.
\end{cases}
\end{equation}
This means that
\begin{equation}\label{4.4}
\begin{cases}
 \textnormal X_t^\epsilon-X_t^z=\epsilon[\e^{tA}\tilde{x}+\int_0^t\e^{(t-s)A}B(\tilde{y}+\tilde{h}(s)-\tilde{h}(0))\d s],\\
  \textnormal Y_t^\epsilon-Y_t^z=\epsilon[\tilde{y}+(\tilde{h}(t)-\tilde{h}(0))].
 \end{cases}
\end{equation}

To construct a coupling $((X_t^z, Y_t^z),(X_t^\epsilon,Y_t^\epsilon))$ by change of measure for them such that
$(X_T^z, Y_T^z)=(X_T^\epsilon,Y_T^\epsilon)$, we take $\tilde{h}$ as follows
\beqlb\label{4.5}
\tilde{h}(t)=\alpha_1(t)\tilde{y}-\alpha_2(t)B^*\e^{(T-t)A^*}U_T^{-1}\e^{TA}\tilde{x}
-\alpha_2(t)B^*\e^{(T-t)A^*}U_T^{-1}\int_0^T\alpha_1(s)\e^{(T-s)A}B\tilde{y}\d s,
\eeqlb
where
$U_T=\int_0^T\alpha_2(t)\e^{(T-t)A}BB^*\e^{(T-t)A^*}\d t,\alpha_i\in C^1([0,T]),i=1,2$ satisfy
$$\alpha_1(0)=1, \alpha_1(T)=0,\alpha_2(0)=\alpha_2(T)=0,\alpha_2(t)>0, \forall t\in(0,T).$$

Note that $U_T$ is invertible due to \cite{Saloff-Coste94} and moreover $\tilde{h}(0)=\tilde{y}, \tilde{h}(T)=0$.

Next, we rewrite the equation \eqref{4.2} as
\begin{equation}\label{4.6}
\begin{cases}
 \textnormal\d X_t^\epsilon=(AX_t^\epsilon+BY_t^\epsilon)\d t,\ \ \ \ \ \ \ \ \ \ \ \ X_0^\epsilon=x+\epsilon\tilde{x}\in\R^{d_1},\\
 \textnormal\d Y_t^\epsilon=\tilde{b}_t(X_t^\epsilon,Y_t^\epsilon)\d t+\sigma(t)\d\bar{B}_t^H, \ \ Y_0^\epsilon=y+\epsilon\tilde{y}\in\R^{d_2}.
\end{cases}
\end{equation}
where
\begin{eqnarray*}
\bar{B}_t^H&:=&B_t^H+\int_0^t\left\{\vartheta(r)\left[\tilde{b}_r(X_r^z,Y_r^z)-\tilde{b}_r(X_r^\epsilon,Y_r^\epsilon)+\epsilon\tilde{h}'(r)\right]\right\}\d r\\
&=&\int_0^tK_H(t,s)\left(\d W_s+\eta_\epsilon(s)\d s\right)
\end{eqnarray*}
and
\beqnn
\eta_\epsilon(s):=
K_H^{-1}\left(\int_0^\cdot\vartheta(r)\left(\tilde{b}_r(X_r^z,Y_r^z)-\tilde{b}_r(X_r^\epsilon,Y_r^\epsilon)+\epsilon\tilde{h}'(r)\right)\d r\right)(s)
\eeqnn
Set
\beqnn
R_\epsilon=\exp\left[-\int_0^T\langle\eta_\epsilon(s),\d W_s\rangle-\frac{1}{2}\int_0^T|\eta_\epsilon(s)|^2\d s\right].
\eeqnn
In the sprit of the proof of Theorem \ref{Th3.2}, we deduce
\beqnn
\int_0^T|\eta_\epsilon(s)|^2\d s\leq C\epsilon^2\left(1+\|B^H\|_\lambda^{2\gamma}\right),
\eeqnn
where and in what follows, $C$ denotes a generic constant.
Consequently, we have
\beqnn
\E\exp\left[\frac{1}{2}\int_0^T|\eta_\epsilon(s)|^2\d s\right]\leq C\E\exp\left[C\epsilon^2\|B^H\|_\lambda^{2\gamma}\right].
\eeqnn
Then by the Fernique theorem, for $\epsilon$ small enough there holds $\E\exp[\frac{1}{2}\int_0^T|\eta_\epsilon(s)|^2\d s]<\infty$.
So, due to the Girsanov theorem for the fractional Brownian motion
(see, e.g., \cite[Theorem 4.9]{Decreusefond&Ustunel98a} or \cite[Theorem 2]{Nualart&Ouknine02b}),
$\{\bar{B}_t^H\}_{t\in[0,T]}$ is a $d$-dimensional fractional Brownian motion under the probability $R_\epsilon\P$.
Therefore, in view of  $(X_T^z, Y_T^z)=(X_T^\epsilon,Y_T^\epsilon)$,  it follows that
\beqnn
P_T f(z+\epsilon\tilde{z})=\E_{R_\epsilon\P}f(X_T^\epsilon,Y_T^\epsilon)=\E\left[R_\epsilon f(X_T^z,Y_T^z)\right].
\eeqnn
Then we arrive at
\beqnn
\frac{P_T f(z+\epsilon\tilde{z})-P_T f(z)}{\epsilon}=\E\left[f(X_T^z,Y_T^z)\frac{R_\epsilon-1}{\epsilon}\right].
\eeqnn
Observe that,
\beqnn
\lim\limits_{\epsilon\rightarrow0}\E\frac{R_\epsilon-1}{\epsilon}
=\lim\limits_{\epsilon\rightarrow0}\E\frac{-\int_0^T\langle\eta_\epsilon(s),\d W_s\rangle-\frac{1}{2}\int_0^T|\eta_\epsilon(s)|^2\d s}{\epsilon}
=\lim\limits_{\epsilon\rightarrow0}\E\frac{-\int_0^T\langle\eta_\epsilon(s),\d W_s\rangle}{\epsilon},
\eeqnn
and moreover by \eqref{2.1},
\beqnn
&&-\int_0^t\langle\eta_\epsilon(s),\d W_s\rangle\\
&=&\frac{H-\frac{1}{2}}{\Gamma(\frac{3}{2}-H)}
\Bigg\{\frac{1}{H-\frac{1}{2}}\int_0^t\left\langle s^{\frac{1}{2}-H}
\vartheta(s)\left(\tilde{b}_s(X_s^\epsilon,Y_s^\epsilon)-\tilde{b}_s(X_s^z,Y_s^z)-\epsilon\tilde{h}'(s)\right),\d W_s\right\rangle\\
&+&\int_0^t\left\langle s^{H-\frac{1}{2}}
\vartheta(s)\left(\tilde{b}_s(X_s^\epsilon,Y_s^\epsilon)-\tilde{b}_s(X_s^z,Y_s^z)-\epsilon\tilde{h}'(s)\right)
\int_0^s\frac{s^{\frac{1}{2}-H}-r^{\frac{1}{2}-H}}{(s-r)^{\frac{1}{2}+H}}\d r,\d W_s\right\rangle\\
&+&
\int_0^t\left\langle\left(\tilde{b}_s(X_s^\epsilon,Y_s^\epsilon)-\tilde{b}_s(X_s^z,Y_s^z)-\epsilon\tilde{h}'(s)\right)
\int_0^s\left(\frac{s}{r}\right)^{H-\frac{1}{2}}
\frac{\vartheta(s)-\vartheta(r)}{(s-r)^{\frac{1}{2}+H}}\d r,\d W_s\right\rangle\\
&+&\int_0^t\left\langle\int_0^s\left(\frac{s}{r}\right)^{H-\frac{1}{2}}
\vartheta(r)
\frac{\left(\tilde{b}_s(X_s^\epsilon,Y_s^\epsilon)-\tilde{b}_s(X_s^z,Y_s^z)\right)
-\left(\tilde{b}_r(X_r^\epsilon,Y_r^\epsilon)-\tilde{b}_r(X_r^z,Y_r^z)\right)}{(s-r)^{\frac{1}{2}+H}}\d r,\d W_s\right\rangle\\
&-&\epsilon\int_0^t\left\langle\int_0^s\left(\frac{s}{r}\right)^{H-\frac{1}{2}}
\vartheta(r)
\frac{\tilde{h}'(s)-\tilde{h}'(r)}{(s-r)^{\frac{1}{2}+H}}\d r,\d W_s\right\rangle
\Bigg\}\\
&=:&\frac{H-\frac{1}{2}}{\Gamma(\frac{3}{2}-H)}[\chi_1(t)+\chi_2(t)+\chi_3(t)+\chi_4(t)+\chi_5(t)].
\\
\eeqnn
Using the B.D.G. inequality and the dominated convergence theorem along with Lemma \ref{Le3.1},
we obtain that as $\epsilon$ goes to zero, $\chi_i(T)/\epsilon, 1\leq i\leq4$, converge to
\beqnn
\frac{1}{H-\frac{1}{2}}\int_0^T\left\langle s^{\frac{1}{2}-H}
\vartheta(s)\left(\nabla\tilde{b}_s(X_s^z,Y_s^z)(\kappa(s),\tilde{h}(s))-\tilde{h}'(s)\right),\d W_s\right\rangle,
\eeqnn
\beqnn
\int_0^T\left\langle s^{H-\frac{1}{2}}
\vartheta(s)\left(\nabla\tilde{b}_s(X_s^z,Y_s^z)(\kappa(s),\tilde{h}(s))-\tilde{h}'(s)\right)
\int_0^s\frac{s^{\frac{1}{2}-H}-r^{\frac{1}{2}-H}}{(s-r)^{\frac{1}{2}+H}}\d r,\d W_s\right\rangle,
\eeqnn
\beqnn
\int_0^T\left\langle\left(\nabla\tilde{b}_s(X_s^z,Y_s^z)(\kappa(s),\tilde{h}(s))-\tilde{h}'(s)\right)
\int_0^s\left(\frac{s}{r}\right)^{H-\frac{1}{2}}
\frac{\vartheta(s)-\vartheta(r)}{(s-r)^{\frac{1}{2}+H}}\d r,\d W_s\right\rangle
\eeqnn
and
\beqnn
\int_0^T\left\langle\int_0^s\left(\frac{s}{r}\right)^{H-\frac{1}{2}}
\vartheta(r)
\frac{\nabla\tilde{b}_s(X_s^z,Y_s^z)(\kappa(s),\tilde{h}(s))
-\nabla\tilde{b}_r(X_r^z,Y_r^z)(\kappa(r),\tilde{h}(r))}{(s-r)^{\frac{1}{2}+H}}\d r,\d W_s\right\rangle
\eeqnn
in $L^1(\P)$, respectively.
The convergence of $\chi_5(T)/\epsilon$ is clear.
Then the result follows.
\fin

\bremark\label{R4.2'}
By \eqref{2.1}, it is not difficult to verify that the right side of \eqref{Th4.1-02} equals to the following expression
\beqnn
\int_0^T\left\langle K_H^{-1}\left(\int_0^\cdot\vartheta(r)\left(\nabla\tilde{b}_r(X_r^z,Y_r^z)\left(\e^{rA}\tilde{x}+\int_0^r\e^{(r-s)A}B\tilde{h}(t)\d t,\tilde{h}(r)\right)-\tilde{h}'(r)\right)\d r\right)(s),\d W_s\right\rangle,
\eeqnn
which is the same as Theorem \ref{Th3.2}.
\eremark

\bremark\label{R4.2}
It is surprise to find that the choosing of $\tilde{h}$ is the same as $\tilde{g}$ in the previous section.
This means that there is some relation between the Malliavin calculus and the coupling argument in dealing with the derivative formula,
which maybe due to the small time asymptotic behaviors investigated by them.
\eremark

\section{Some Applications: (log) Harnack Inequalities and The Hyperbounded Property}

\setcounter{equation}{0}

In this section, we give some applications of the derivative formulas.
To derive the explicit Harnack type inequalities and the hyperbounded property from Theorem \ref{Th4.1-01},
we let
\beqlb\label{4.7}
\tilde{h}(t)=\alpha_1(t)\tilde{y}-\alpha_2(t)B^*\e^{(T-t)A^*}U_T^{-1}\e^{TA}\tilde{x}
\eeqlb
for
\beqnn
\alpha_1(0)=1,\ \alpha_1(T)=0,\ \int_0^T\alpha_1(t)\e^{(T-t)A}B\d t=0,
\eeqnn
and
\beqnn
\alpha_2(0)=\alpha_2(T)=0, \alpha_2(t)>0, \forall t\in(0,T).
\eeqnn
It is not hard to verify that under $\tilde{h}$ there also holds $(X_T^z, Y_T^z)=(X_T^\epsilon,Y_T^\epsilon)$.
Besides, we assume that $A^{n_0}=0$ for some $n_0\in\mathbb{Z}^+$.
As a consequence, we may choose $\alpha_i,i=1,2$ as follows
$$\alpha_1(t)=\sum_{i=1}^{n_0+1}a_i\frac{(T-t)^i}{T^i},\ \ \alpha_2(t)=\frac{t(T-t)}{T^2},$$
where the coefficients $a_i,1\leq i\leq n_0+1$ satisfy
\begin{equation}\label{1.2}\nonumber
\left\{
\begin{array}{ll}
\sum_{i=1}^{n_0+1}a_i=1,\\
\sum_{i=1}^{n_0+1}a_i\frac{1}{i+j+1}=0,\ 1\leq j\leq n_0+1.
\end{array} \right.
\end{equation}

\bremark\label{R4.1}
For \eqref{4.7}, the choosing of $\alpha_1$ is not unique.
For instance, we may take $\alpha_1$ another form
$$\alpha_1(t)=1+\sum_{i=1}^{n_0+1}\tilde{a}_i\frac{t^i}{T^i},$$
where the coefficient $\tilde{a}_i,1\leq i\leq n_0+1$ satisfy
\begin{equation}\label{1.2}\nonumber
\left\{
\begin{array}{ll}
1+\sum_{i=1}^{n_0+1}\tilde{a}_i=0,\\
1+\sum_{i=1}^{n_0+1}\tilde{a}_i\frac{1}{C_{i+j+1}^{j+1}}=0,\ 1\leq j\leq n_0+1.
\end{array} \right.
\end{equation}
\eremark

Besides, we make use of the following Kalman rank condition (see also Remark \ref{addRemark(Th3.2)}): (A) \\
There exists an integer number $k_0\in[0,d_1-1]$ such that Rank$[B,AB,\cdots,A^{k_0}B]=d_1$.

Using the above assumptions, we get the following Harnack type inequalities and the hyperbounded property.
\btheorem\label{Th4.1-002}
Assume (H1)(ii), (H2) and (A).
Then,
\begin{itemize}
\item[(1)] for any nonnegative $f\in\mathcal{B}_b(\R^{d_1+d_2})$,
\beqnn
(P_T f(z))^p&\leq&P_Tf^p(z+\tilde{z})
\exp\Bigg[\frac{p}{p-1}\Bigg(a(T)|\tilde{x}|^2+\tilde{a}(T)|\tilde{y}|^2\\
&&+\left(1+\left(1\vee\left(\frac{p^2}{(p-1)^2}\left(b(T)|\tilde{x}|^2+\tilde{b}(T)|\tilde{y}|^2\right)\right)\right)^{\frac{\gamma}{1-\gamma}}
\right)\left(b(T)|\tilde{x}|^2+\tilde{b}(T)|\tilde{y}|^2\right)\Bigg)
\Bigg];
\eeqnn

\par

\item[(2)] for any positive $f\in\mathcal{B}_b(\R^{d_1+d_2})$,
\beqnn
P_T\log f(z)&\leq&\log P_T f(z+\tilde{z})
+\Bigg[a(T)|\tilde{x}|^2+\tilde{a}(T)|\tilde{y}|^2\\
&&+\left(1+\left(1\vee\left(\frac{p^2}{(p-1)^2}\left(b(T)|\tilde{x}|^2+\tilde{b}(T)|\tilde{y}|^2\right)\right)\right)^{\frac{\gamma}{1-\gamma}}
\right)\left(b(T)|\tilde{x}|^2+\tilde{b}(T)|\tilde{y}|^2\right)\Bigg];
\eeqnn

\par

\item[(3)] if $\mu$ is a probability measure on $\R^{d_1+d_2}$
such that for some $\tilde{K}(>0), \mu(P_T f)\leq\tilde{K}\mu(f), \forall f\in\mathcal{B}_b^+(\R^{d_1+d_2})$,
there holds
\beqnn
\|P_T\|_{p\rightarrow\upsilon p}^{\upsilon p}\leq\tilde{K}^\upsilon\int_{\R^{d_1+d_2}}\frac{\mu(\d z)}{\left(\int_{\R^{d_1+d_2}}\e^{-\Phi(z,\hat{z})}\mu(\d\hat{z})\right)^\upsilon},\ \forall\upsilon>1,
\eeqnn
with
\beqnn
\Phi(z,\hat{z})&=&\frac{p}{p-1}\Bigg[a(T)|(\hat{z}-z)_{d_1}|^2+\tilde{a}(T)| (\hat{z}-z)_{d_2}|^2\\
&&+\left(1+\left(1\vee\left(\frac{p^2}{(p-1)^2}\left(b(T)|(\hat{z}-z)_{d_1}|^2+\tilde{b}(T)|(\hat{z}-z)_{d_2}|^2\right)\right)\right)^{\frac{\gamma}{1-\gamma}}
\right)\\
&&~~~\times\left(b(T)| (\hat{z}-z)_{d_1}|^2+\tilde{b}(T)| (\hat{z}-z)_{d_2}|^2\right)\Bigg],
\eeqnn
\end{itemize}
where $z,\tilde{z}\in\R^{d_1+d_2}$ and $a(T),\tilde{a}(T),b(T),\tilde{b}(T)$ are defined below,
$(\hat{z}-z)_{d_1}$ and $(\hat{z}-z)_{d_2}$ stand for the first and the second components of $\hat{z}-z$, respectively.
\etheorem

\emph{Proof.}
By \cite[Theorem 4.2]{Wang&Zhang13}, the condition (A) and the expression of $\alpha_2$ yield that
$$\|U_T^{-1}\|\leq C\frac{T}{(T\wedge1)^{2(k_0+1)}}.$$
Without lost of generality, we assume $T\leq1$.
Then there holds $\|U_T^{-1}\|\leq C\frac{1}{T^{2k_0+1}}.$\\
Note that, for each $t,s\in[0,T]$ and $i=1,2$,
\beqnn
|\alpha_i(t)|\leq C,\ |\alpha_i'(t)|\leq C\frac{1}{T},\ |\alpha_i'(t)-\alpha_i'(s)|\leq C\frac{1}{T^2}|t-s|.
\eeqnn
Consequently, it is easy to see that
\beqnn
|(\kappa(t),\tilde{h}(t))|\leq C\left[\left(1+\frac{1}{T^{2k_0+1}}\right)|\tilde{x}|+|\tilde{y}|\right],
\eeqnn
\beqnn
|\tilde{h}'(t)|\leq C\left[\left(\frac{1}{T^{2k_0+1}}+\frac{1}{T^{2(k_0+1)}}\right)|\tilde{x}|+\frac{1}{T}|\tilde{y}|\right]
\eeqnn
and
\beqnn\tilde{}
|\tilde{h}'(t)-\tilde{h}'(s)|\leq\left[\frac{1}{T^{2k_0+1}}\left(1+\frac{1}{T}+\frac{1}{T^2}\right)|\tilde{x}|+\frac{1}{T^2}|\tilde{y}|\right]|t-s|
\eeqnn
hold for all $t,s\in[0,T]$.
Combining these with the expression of $M_T$, we deduce that
\beqlb\label{Pf(Th4.1-02)-1}
\langle M\rangle_T\leq
\left(a(T)|\tilde{x}|^2+\tilde{a}(T)|\tilde{y}|^2\right)+\left(b(T)|\tilde{x}|^2+\tilde{b}(T)|\tilde{y}|^2\right)\|B^H\|_\lambda^{2\gamma},
\eeqlb
where
\beqnn
a(T)&=&CT^{2-2H}\Bigg\{1+T^{2\delta}+T^{2\gamma}+T^{2\varrho}+T^2\\
&&+\frac{1}{T^{4k_0}}\left[T^2+1+\frac{1}{T^{2(1-\delta)}}
+\frac{1}{T^{2(1-\gamma)}}
+\frac{1}{T^{2(1-\varrho)}}
+\frac{1}{T^{2(2-\delta)}}
+\frac{1}{T^2}+\frac{1}{T^4}\right]\Bigg\},
\eeqnn
\beqnn
\tilde{a}(T)=CT^{2-2H}\left[1+T^{2\delta}+T^{2\gamma}+T^{2\varrho}+T^2+T^4+\frac{1}{T^{2(1-\delta)}}+\frac{1}{T^2}\right],
\eeqnn
and
\beqnn
b(T)=CT^{2(\lambda\gamma-H+1)}\left(1+\frac{1}{T^{2(2k_0+1)}}\right),\ \tilde{b}(T)=CT^{2(\lambda\gamma-H+1)}.
\eeqnn
On the other hand, by Theorem \ref{Th4.1-01} and the Young inequality (see, e.g., \cite[Lemma 2.4]{Arnaudon&Thalmaier&Wang09a}),
we obtain that, for all $\theta>0$,
\beqlb\label{Pf(Th4.1-02)-2}
&&|\nabla_{\tilde{z}}{P_T f(z)}|-\theta[P_T(f\log f)(z)-(P_T f)(z)(\log P_T f)(z)]\cr
&\leq&\theta\log\E\exp\left[\frac{1}{\theta}M_T\right]\cdot P_T f(z)
\leq\frac{\theta}{2}\log\E\exp\left[\frac{2}{\theta^2}\langle M\rangle_T\right]\cdot P_T f(z).
\eeqlb
In the sprit of \cite[Lemma 3.3]{Fan15a}, by \eqref{Pf(Th4.1-02)-1} we conclude that
\beqnn
&&\E\exp\left[\frac{2}{\theta^2}\langle M\rangle_T\right]\\
&\leq&\exp\Bigg[\frac{1}{\theta^2}\Bigg(a(T)|\tilde{x}|^2+\tilde{a}(T)|\tilde{y}|^2\\
&&~~~~~~+\left(1+\left(1\vee\left(\frac{p^2}{(p-1)^2}\left(b(T)|\tilde{x}|^2+\tilde{b}(T)|\tilde{y}|^2\right)\right)\right)^{\frac{\gamma}{1-\gamma}}
\right)\left(b(T)|\tilde{x}|^2+\tilde{b}(T)|\tilde{y}|^2\right)
\Bigg)\Bigg].
\eeqnn
This, together with \eqref{Pf(Th4.1-02)-2}, implies that
\beqlb\label{Pf(Th4.1-02)-3}
&&|\nabla_{\tilde{z}}{P_T f(z)}|-\theta[P_T(f\log f)(z)-(P_T f)(z)(\log P_T f)(z)]\cr
&\leq&
\Bigg[a(T)|\tilde{x}|^2+\tilde{a}(T)|\tilde{y}|^2\cr
&&+\left(1+\left(1\vee\left(\frac{p^2}{(p-1)^2}\left(b(T)|\tilde{x}|^2+\tilde{b}(T)|\tilde{y}|^2\right)\right)\right)^{\frac{\gamma}{1-\gamma}}
\right)\left(b(T)|\tilde{x}|^2+\tilde{b}(T)|\tilde{y}|^2\right)
\Bigg]\frac{1}{\theta}P_T f(z).\nonumber\\
\cr
\eeqlb
Hence, by \cite[Proposition 1.3.1]{Wang13a} we show the first assertion.
While the second assertion then follows from \cite[Corollary 1.3.4]{Wang13a} by noting that $\R^{d_1+d_2}$ is a length space.
Finally, for $f\in\mathcal{B}_b^+(\R^{d_1+d_2})$ with $\int_{\R^{d_1+d_2}}f^p\d\mu=1$, by the inequality in (1) we get
\beqnn
(P_T f(z))^p\e^{-\Phi(z,\hat{z})}\leq P_Tf^p(\hat{z}),
\eeqnn
where
\beqnn
\Phi(z,\hat{z})&=&\frac{p}{p-1}\Bigg[a(T)|(\hat{z}-z)_{d_1}|^2+\tilde{a}(T)| (\hat{z}-z)_{d_2}|^2\\
&&+\left(1+\left(1\vee\left(\frac{p^2}{(p-1)^2}\left(b(T)|(\hat{z}-z)_{d_1}|^2+\tilde{b}(T)|(\hat{z}-z)_{d_2}|^2\right)\right)\right)^{\frac{\gamma}{1-\gamma}}
\right)\\
&&~~~\times\left(b(T)| (\hat{z}-z)_{d_1}|^2+\tilde{b}(T)| (\hat{z}-z)_{d_2}|^2\right)\Bigg].
\eeqnn
Then integrating both sides with respect to $\mu(\d\hat{z})$ implies
\beqnn
(P_T f(z))^p\int_{\R^{d_1+d_2}}\e^{-\Phi(z,\hat{z})}\mu(\d\hat{z})\leq\mu(P_Tf^p)\leq\widetilde{K}\mu(f^p)=\tilde{K}.
\eeqnn
Consequently, for every $\upsilon>1$, we have
\beqnn
(P_T f(z))^{\upsilon p}\leq\frac{\tilde{K}^\upsilon}{\left(\int_{\R^{d_1+d_2}}\e^{-\Phi(z,\hat{z})}\mu(\d\hat{z})\right)^\upsilon}.
\eeqnn
Then, we conclude that
\beqnn
\int_{\R^{d_1+d_2}}(P_T f(z))^{\upsilon p}\mu(\d z)\leq
\tilde{K}^\upsilon\int_{\R^{d_1+d_2}}\frac{\mu(\d z)}{\left(\int_{\R^{d_1+d_2}}\e^{-\Phi(z,\hat{z})}\mu(\d\hat{z})\right)^\upsilon}.
\eeqnn
\fin

\bremark\label{Re5.1}
Going back to the above proof,
we also obtain the gradient-entropy inequality, that is \eqref{Pf(Th4.1-02)-3}.
Besides, as a direct application of the Harnack type inequalities,
according to \cite[Proposition 4.1]{Prato&Rockner&Wang09a} we get the strong Feller property of $P_T$.
That is,
the expression $\lim_{|\tilde{z}-z|\rightarrow0}P_T f(\tilde{z})=P_T f(z)$ holds for each $f\in\mathcal{B}_b(\R^{d_1+d_2})$ and $z\in\R^{d_1+d_2}$.
For the applications of these inequalities to invariant probability measure, Entropy-cost inequalities and heat kernel estimates,
one can refer to the monograph \cite{Wang13a}.
\eremark

\bremark\label{Re5.2}
As in \cite{Zhang10a,Wang&Zhang13},
we may generalize the above results to the following stochastic Hamiltonian system driven by fractional Brownian motion:
\begin{equation}\label{6.1}
\begin{cases}
 \textnormal\d X_t=\nabla^2\mathcal {H}(X_t,Y_t)\d t,\ \ \ \ \ \ \ \ \ \ \ \ \ \ \ \ \ \ \ ~X_0=x\in\R^{d},\nonumber\\
 \textnormal\d Y_t=-\nabla^1\mathcal {H}(X_t,Y_t)\d t+\sigma(t)\d B_t^H, \ ~~~Y_0=y\in\R^{d},
\end{cases}
\end{equation}
where $\mathcal {H}$ is a Hamiltonian function.
\eremark

\textbf{Acknowledgement}
The author would like to thank the referees for detailed comments and valuable suggestions that helped to improve the paper significantly.


\begin{thebibliography}{99}


\bibitem{Alos&Mazet&Nualart01a} E. Al\`{o}s, O. Mazet and D. Nualart,  {\it Stochastic calculus with respect to Gaussian processes}, Ann. Probab.
29(2001), 766--801.

\bibitem{Arnaudon&Thalmaier&Wang09a} M. Arnaudon, A. Thalmaier and F. Y. Wang, {\it Gradient estimates and Harnack inequalities on non-compact Riemannian
  manifolds}, Stochastic Process. Appl. 119(2009), 3653--3670.

\bibitem{Xiao11a} A. Ayache, N. R. Shieh and Y. M. Xiao, {\it Multiparameter multifractional Brownian motion:
local nondeterminism and joint continuity of the local times}, Ann. Inst. H. Poincar\'{e}
Probab. Statist. 47(2011), 1029--1054.



\bibitem{Bao&Wang&Yuan13b} J. Bao, F. Y. Wang and C. Yuan,  {\it Derivative formula and Harnack inequality for degenerate functional SDEs},
Stochastics and Dynamics 13(2013), 1--22.

\bibitem{Bao&Wang&Yuan13c} J. Bao, F. Y. Wang and C. Yuan, {\it Bismut formulae and applications for functional SPDEs},
Bull. Sci. Math. 137(2013), 509--522.

\bibitem{Baudoin&Ouyang11} F. Baudoin and C. Ouyang, {\it Small-time kernel expansion for solutions of stochastic differential equations driven by fractional Brownian motions}, Stochastic Process Appl. 121(2011), 759--792.

\bibitem{Baudoin&Ouyang&Tindel11a} F. Baudoin, C. Ouyang and S. Tindel, {\it Upper bounds for the density of solutions of stochastic differential equations driven by fractional Brownian motions}, Ann. Inst. H. Poincar¨¦ Probab. Statist. 50(2014), 111--135.

\bibitem{Biagini&Hu08a} F. Biagini, Y. Hu, B. $\emptyset$ksendal and  T. Zhang, {\it Stochastic Calculus for Fractional Brownian Motion and Applications}, Springer-Verlag, London, 2008.

\bibitem{Bismut84} J. M. Bismut, \textit{Large Deviation and The Malliavin Calculus}, Birkh\"{a}user, Boston, MA,1984.


\bibitem{Coutin&Qian02a} L. Coutin and Z. Qian, {\it Stochastic analysis, rough path analysis and fractional Brownian motions},
   Probab. Theory Related Fields 122(2002), 108--140.


\bibitem{Decreusefond&Ustunel98a} L. Decreusefond and A. S. \"{U}st\"{u}nel, {\it Stochastic analysis of the fractional Brownian motion}, Potential Anal.
    10(1998), 177--214.

\bibitem{Dong&Xie10} Z. Dong and Y. Xie,  {\it Ergodicity of linear SPDE driven by L\'{e}vy noise},
J. Syst. Sci. Complex 23(2010), 137--152.

\bibitem{Driver97} B. Driver, \textit{Integration by parts for heat kernel measures revisited}, J. Math. Pures Appl. 76(1997), 703--737.


\bibitem{Elworthy&Li94} K. D. Elworthy and X. M. Li, {\it Formulae for the derivatives of heat semigroups}, J. Funct. Anal. 125(1994), 252--286.

\bibitem{Fan13b} X. L. Fan, {\it Harnack inequality and derivative formula for SDE driven by fractional Brownian motion},
  Science in China-Mathematics 561(2013), 515--524.

\bibitem{Fan15a} X. L. Fan, {\it Integration by parts formula and applications for SDEs driven by fractional Brownian motions}, Stochastic Analysis and Applications 33(2015), 199--212.

\bibitem{Fan15b} X. L. Fan, {\it Logarithmic Sobolev inequalities for fractional diffusion},  Statist. Probab. Lett. 106(2015), 165--172.


\bibitem{Fan13a} X. L. Fan and Y. Ren, {\it Bismut formulae and applications for stochastic (functional) differential equations driven by fractional Brownian motions}, Stochastics and Dynamics. 7(2017), 19 pages.



\bibitem{Fang&Li&Luo11} S. Fang, H. Li and D. Luo, {\it Heat semi-group and generalized flows on complete Riemannian manifolds},
 Bull. Sci. Math., 135(2011), 565--600.


\bibitem{Guillin&Wang12} A. Guillin and F. Y. Wang,  {\it Degenerate Fokker-Planck equations: Bismut formula, gradient estimate and Harnack inequality},
J. Differential Equations 253(2012), 20--40.


\bibitem{Hairer05} M. Hairer, {\it Ergodicity of stochastic differential equations driven by fractional Brownian motion}, Ann. Probab. 33(2005), 703--758.

\bibitem{Hairer&Pillai11} M. Hairer and N. S. Pillai, {\it Ergodicity of hypoelliptic SDEs driven by fractional Brownian motion}, Ann. Inst. H. Poincar\'{e} Probab. Statist. 47(2011), 601--628.



\bibitem{Jaramillo&Nualart17a} A. Jaramillo and D. Nualart, {\it Asymptotic properties of the derivative of self-intersection local time of fractional Brownian motion},  Stochastic Process. Appl. 127(2017), 669--700.


\bibitem{Nikiforov&Uvarov88} A. F. Nikiforov and V. B. Uvarov, {\it Special Functions of Mathematical Physics}, Birkh\"{a}user, Boston, 1988.

\bibitem{Nourdin&Simon06} I. Nourdin and T. Simon, {\it On the absolute continuity of one-dimensional SDEs driven by a fractional Brownian motion}, Statist. Probab. Lett. 76(2006), 907--912.

\bibitem{Nualart06a} D. Nualart, {\it The Malliavin Calculus and Related Topics, Second edition},
Springer-Verlag, Berlin, 2006.

\bibitem{Nualart&Ouknine02b} D. Nualart and Y. Ouknine,  {\it Regularization of differential equations by fractional noise},  Stochastic Process. Appl.
    102(2002), 103--116.

\bibitem{Nualart&Rascanu02a} D. Nualart and A. R\u{a}\c{s}canu, {\it Differential equations driven by fractional Brownian motion}, Collect. Math. 53(2002), 55--81.

\bibitem{Nualart&Saussereau09} D. Nualart and B. Saussereau, {\it Malliavin calculus for stochastic differential equations driven by a fractional Brownian motion}, Stochastic Process. Appl. 119(2009), 391--409.



\bibitem{Prato&Rockner&Wang09a} G. Da Prato, M. R\"{o}ckner and F.Y. Wang, {\it Singular stochastic equations on Hilberts space: Harnack inequalities for their transition semigroups}, J. Funct. Anal. 257(2009), 992--1017.

\bibitem{Priola06} E. Priola,  {\it Formulae for the derivatives of degenerate diffusion semigroups},
J. Evol. Equ. 6(2006), 577--600.


\bibitem{Saloff-Coste94} L. Saloff-Coste, {\it Convergence to equilibrium and logarithmic Sobolev constant on manifolds with Ricci
curvature bounded below}, Colloq. Math. 67(1994), 109--121.

\bibitem{Samko&Kilbas&Marichev} S. G. Samko, A. A. Kilbas and O. I. Marichev,  {\it Fractional Integrals and Derivatives, Theory and Applications}, Gordon and Breach Science Publishers, Yvendon, 1993.

\bibitem{Saussereau12} B. Saussereau, {\it Transportation inequalities for stochastic differential equations driven
by a fractional Brownian motion}, Bernoulli 18(2012), 1--23.


\bibitem{Wang11a} F. Y. Wang, {\it Derivative formula and Harnack inequality for SDEs driven by L\'{e}vy processes},
Stochastic Analysis and Applications 32(2014), 30--49.


\bibitem{Wang13a} F. Y. Wang, {\it Harnack Inequalities for Stochastic Partial Differential Equations}, Springer, 2013.

\bibitem{Wang&Xu10a} F. Y. Wang and L. Xu, {\it Derivative formulae and applications for hyperdissipative stochastic Navier-Stokes/Burgers equations},
Infin. Dimens. Anal. Quantum Probab. Relat. Top. 15(2012), 1--19.

\bibitem{Wang&Zhang13} F. Y. Wang and X. C. Zhang, {\it Derivative formula and applications for degenerate diffusion semigroups},
J. Math. Pures Appl. 99(2013), 726--740.


\bibitem{Yan16a} L. T. Yan, {\it The fractional derivative for fractional Brownian local time},
Math. Z. 283(2016), 437--468.


\bibitem{Young36a} L. C. Young,  {\it An inequality of the H\"{o}lder type connected with Stieltjes integration}, Acta Math. 67(1936), 251--282.


\bibitem{Zahle98a} M. Z\"{a}hle, {\it Integration with respect to fractal functions and stochastic calculus I}, Probab. Theory Related Fields 111(1998),
333--374.

\bibitem{ZhangSQ13a} S. Q. Zhang, {\it Harnack inequality for semilinear SPDEs with multiplicative noise}, Statist. Probab. Lett. 83(2013), 1184--1192.

\bibitem{Zhang10a} X. C. Zhang, {\it Stochastic flows and Bismut formulas for stochastic Hamiltonian systems}, Stochastic Process. Appl. 120(2010), 1929--1949.


\end{thebibliography}
\end{document}